\documentclass[12pt,draft]{amsart}
\usepackage{amssymb, amsmath}
\usepackage{url}

\newtheorem{theorem}{Theorem}[section]

\theoremstyle{remark}

\theoremstyle{definition}

\numberwithin{equation}{section}
\numberwithin{theorem}{section}

\newcommand{\intinf}{\int^\infty_{-\infty}}

\newcommand{\R}{{\mathbb R}}

\newcommand{\fn}{\!:\!}

\providecommand{\norm}[1]{\lVert#1\rVert}

\newcommand{\sqrtpt}{\sqrt{\pi t}}

\begin{document}


\date{Preprint January 2, 2023.}
\title[Heat equation]{Sharp norm estimates for the classical heat equation}
\author{Erik Talvila}
\address{Department of Mathematics \& Statistics\\
University of the Fraser Valley\\
Abbotsford, BC Canada V2S 7M8}
\email{Erik.Talvila@ufv.ca}

\begin{abstract}
Sharp estimates of solutions of the classical heat equation are proved
in $L^p$ norms on the real line.
\end{abstract}

\maketitle

\section{Introduction}\label{sectionintroduction}

In this paper we give sharp estimates of solutions of the classical heat equation on the real line
with initial value data that is in an
$L^p$ space
($1\leq p\leq\infty$).

For $u\fn\R\times(0,\infty)\to\R$ write $u_t(x)=u(x,t)$.

The classical problem of the heat equation on the real line is, given a function
$f\in L^p$ for some $1\leq p\leq\infty$,  find a function
$u\fn \R\times (0,\infty)\to\R$ such that
$u_t\in C^2(\R)$ for each $t>0$, $u(x,\cdot)\in C^1((0,\infty))$ for each $x\in\R$ and
\begin{align}
&\frac{\partial^2u(x,t)}{\partial x^2}-\frac{\partial u(x,t)}{\partial t}=0 
\text{ for each } (x,t)\in\R\times(0,\infty)\label{heatpde}\\
&\lim_{t\to 0^+}\norm{u_t-f}_p=0.\label{Lpic}
\end{align}
If $p=\infty$ then $f$ is also assumed to be continuous.

A solution
is given by 
the convolution $u_t(x)=F\ast\Theta_t(x)=\intinf F(x-y)\Theta_t(y)\,dy$ 
where the Gauss--Weierstrass heat kernel
is $\Theta_t(x)=\exp(-x^2/(4t))/(2\sqrtpt)$.  
For example, see \cite{follandpde}.  Under suitable growth conditions on $u$
the solution is unique.
See \cite{hirschmanwidder} and \cite{widderbook}.  References
\cite{cannon} and \cite{widderbook} contain many results on the classical heat equation, including 
extensive bibliographies.

The heat kernel has the following properties.  Let $t>0$ and let $s\not=0$ such
that $1/s+1/t>0$.  Then
\begin{align}
\Theta_t\ast\Theta_s&=\Theta_{t+s}\label{thetaconvolution}\\
\norm{\Theta_t}_q&=\frac{\alpha_q}{t^{(1-1/q)/2}} \text{ where } \alpha_q=\left\{\begin{array}{cl}
1, & q=1\\
\frac{1}{(2\sqrt{\pi})^{1-1/q} \,q^{1/(2q)}}, & 1<q<\infty\\
\frac{1}{2\sqrt{\pi}}, & q=\infty.
\end{array}
\right.\label{Thetaqnormalpha}
\end{align}

The last of these follows from the probability integral $\intinf e^{-x^2}\,dx=\sqrt{\pi}$.

\begin{theorem}\label{theoremLpestimates}
Let $1\leq p\leq\infty$ and $f\in L^p$.\\
(a) If $p\leq s\leq\infty$ then $f\ast\Theta_t\in L^s$.\\ 
(b)
Let $q,r\in[1,\infty]$
such that $1/p+1/q=1+1/r$.  
There is a constant $K_{p,q}$ such
that $\norm{f\ast\Theta_t}_r\leq K_{p,q}\norm{f}_p \,t^{-(1-1/q)/2}$ for all $t>0$.  The estimate
is sharp in the sense that if $\psi\fn(0,\infty)\to(0,\infty)$ such that $\psi(t)=o(t^{-(1-1/q)/2})$
as $t\to 0^+$ or $t\to\infty$ then there is $G\in L^p$ such that $\norm{G\ast\Theta_t}_r/\psi(t)$ is
not bounded as $t\to 0^+$ or $t\to\infty$.  The constant 
$K_{p,q}=(c_pc_q/c_r)^{1/2}\alpha_q$, where $c_p=p^{1/p}/(p')^{1/p'}$ with $p, p'$ being
conjugate exponents.  It cannot be replaced with any smaller
number.\\
(c) If $1\leq s<p$ then
$f\ast\Theta_t$ need not be in $L^s$.
\end{theorem}

When $r=p$ and $q=1$ the inequality in part (b) reads $\norm{f\ast\Theta_t}_p\leq\norm{f}_p$.
When $r=\infty$ then $p$ and $q$ are conjugates and 
the inequality in part (b) reads $\norm{f\ast\Theta_t}_\infty\leq\norm{f}_pt^{-1/(2p)}$.

The condition for sharpness in Young's inequality is that both functions be Gaussians.
This fact is
exploited in the proof of part (b).
See \cite[p.~99]{liebloss}, \cite{beckner} and \cite{toscani}.
Our proof also uses ideas from
\cite[Theorem~9.2, p.~195]{hirschmanwidder} and \cite[pp.~115-120]{axler}.

The estimates are known, for example \cite[Proposition~3.1]{iwabuchi}, but we have not been able to find a proof in the
literature that they are sharp.

\begin{proof}
(a), (b)  Young's inequality gives
\begin{equation}
\norm{f\ast\Theta_t}_r\leq C_{p,q}\norm{f}_p\norm{\Theta_t}_q =\frac{C_{p,q}\norm{f}_p\alpha_q}{
t^{(1-1/q)/2}},\label{normyounginequality}
\end{equation}
where $\alpha_q$ is given in \eqref{Thetaqnormalpha}.  
The sharp constant, given in \cite[p.~99]{liebloss},
is $C_{p,q}=(c_pc_q/c_r)^{1/2}$ where $c_p=p^{1/p}/(p')^{1/p'}$ with $p, p'$ being
conjugate exponents. Note that $c_1=c_{\infty}=1$.  Also, $0<C_{p,q}\leq 1$.   We then
take $K_{p,q}=C_{p,q}\alpha_q$.

To show the estimate $\norm{u_t}_r=O(t^{-(1-1/q)/2})$ is sharp as $t\to 0^+$ and $t\to\infty$, let
$\psi$ be as in the statement of the theorem.  Fix $p\leq r\leq\infty$.  Define the family of linear operators
$S_t\fn L^p\to L^r$ by $S_t[f](x)=f\ast \Theta_t(x)/\psi(t)$.  The estimate
$\norm{S_t[f]}_r\leq K_{p,q}\norm{f}_pt^{-(1-1/q)/2}/\psi(t)$ shows that, for each $t>0$, $S_t$
is a bounded linear operator.  Let $f_t=\Theta_t$.  Then, from \eqref{thetaconvolution} and 
\eqref{Thetaqnormalpha},
$$
\frac{\norm{S_t[f_t]}_r}{\norm{f_t}_p}=
\frac{\norm{\Theta_t\ast\Theta_t}_r}{\psi(t)\norm{\Theta_t}_p}
=\frac{\norm{\Theta_{2t}}_r}{\psi(t)\norm{\Theta_t}_p}
=\frac{\alpha_r}{\alpha_p 2^{(1-1/r)/2}\psi(t)t^{(1-1/q)/2}}.
$$
This is not bounded in the limit $t\to 0^+$.  Hence, $S_t$ is not uniformly bounded.
By the Uniform Bounded Principle it is not pointwise bounded.  Therefore, there is a function
$f\in L^p$ such that $\norm{f\ast\Theta_t}_r\not=O(\psi(t))$ as $t\to 0^+$.  And, the growth estimate
$\norm{f\ast\Theta_t}_r=O(t^{-(1-1/q)/2)})$ as $t\to 0^+$ is sharp.  Similarly for sharpness
as $t\to\infty$.

Now show the constant $K_{p,q}$ cannot be reduced.
A calculation shows we have equality in \eqref{normyounginequality}
when $f=\Theta^\beta_t$ and
 $\beta$ is given
by the equation
\begin{equation}
\frac{\beta^{1-1/q}}{(\beta +1)^{1-1/r}}
  =  
\frac{c_pc_q}{c_r}\left(\frac{\alpha_p\alpha_q}{\alpha_r}\right)^2
  =  
\left(1-\frac{1}{p}\right)^{1-1/p}\left(1-\frac{1}{q}\right)^{1-1/q}
\left(1-\frac{1}{r}\right)^{-(1-1/r)}.\label{youngequality}
\end{equation}

First consider the case $p\not=1$ and $q\not=1$.
Notice that $1-1/r=(1-1/q)+(1-1/p)>1-1/q$.  Let
$g(x)=x^A(x+1)^{-B}$ with $B>A>0$.
Then $g$ is
strictly increasing on $(0,A/(B-A))$ and strictly decreasing for $x> A/(B-A)$ so there is a unique
maximum for $g$ at $A/(B-A)$. Put $A=1-1/q$ and $B=1-1/r$.  Then
$$
g\left(\frac{A}{B-A}\right) = \frac{\beta^{1-1/q}}{(\beta +1)^{1-1/r}}
= \left(1-\frac{1}{p}\right)^{1-1/p}\left(1-\frac{1}{q}\right)^{1-1/q}
\left(1-\frac{1}{r}\right)^{-(1-1/r)}. 
$$
Hence, \eqref{youngequality} has a unique positive solution for $\beta$ given by
$\beta=(1-1/q)/(1-1/p)$.  

If $p=1$ then $q=r$.  In this case, \eqref{youngequality} reduces to
$(1+1/\beta)^{1-1/q}=1$ and the solution is given in the limit $\beta\to\infty$.  Sharpness
of \eqref{normyounginequality} is then given in this limit.  It can also be seen that taking $f$ to be the Dirac
distribution gives equality.

If $q=1$ then $p=r$.  Now, \eqref{youngequality} reduces to $(\beta+1)^{1-1/p}=1$ and $\beta=0$.  There
is equality in \eqref{normyounginequality} when $f=1$.  This must be done in the limit $\beta\to0^+$.

If $p=q=r=1$ then there is equality in \eqref{normyounginequality} for each $\beta>0$.

Hence, the constant in \eqref{normyounginequality} is sharp.

(c) Suppose $f\geq 0$ and $f$ is decreasing on $[c,\infty)$ for some $c\in\R$.  Let $x>c$.  Then
\begin{eqnarray*}
f\ast\Theta_t(x) & \geq & \int_c^x f(y)\Theta_t(x-y)\,dy
 \geq  f(x)\int_c^x \Theta_t(x-y)\,dy\\
 & = & \frac{f(x)}{\sqrt{\pi}}\int_0^{(x-c)/(2\sqrt{t})} e^{-y^2}\,dy
 \sim  f(x)/2 \quad \text{ as } x \to \infty.
\end{eqnarray*}
Now put $f(x)=1/[x^{1/p}\log^2(x)]$ for $x\geq e$ and $f(x)=0$, otherwise.
For $p=\infty$ replace $x^{1/p}$ by $1$.
\end{proof}


\begin{thebibliography}{99}
\bibitem{axler}
S. Axler, P. Bourdon and W. Ramey, {\it Harmonic function theory},
New York, Springer-Verlag, 2001.
\bibitem{beckner}
W. Beckner,
{\it Inequalities in Fourier analysis},
Ann. of Math. (2) {\bf 102}(1975), 159--182.
\bibitem{cannon}
J.R. Cannon, {\it The one-dimensional heat equation},
Menlo Park, Addison--Wesley, 1984.
\bibitem{follandpde}
G.B. Folland, {\it Introduction to partial differential equations},
Princeton, Princeton University Press, 1995.
\bibitem{hirschmanwidder}
I.I. Hirschman and D.V. Widder, {\it The convolution transform},
Princeton, Princeton University Press, 1955.
\bibitem{iwabuchi}
T. Iwabuchi, T. Matsuyama and K. Taniguchi,
{\it Boundedness of spectral multipliers for Schr\"odinger operators on open sets},
Rev. Mat. Iberoam. {\bf 34}(2018), 1277--1322.
\bibitem{liebloss}
E.H. Lieb and M. Loss, {\it Analysis}, Providence, American Mathematical
Society, 2001.
\bibitem{toscani}
G. Toscani,
{\it Heat equation and the sharp Young’s inequality},
arXiv:1204.2086 (2012).
\bibitem{widderbook}
D.V. Widder, {\it The heat equation}, New York, Academic Press, 1975.
\end{thebibliography}
\end{document}